\documentstyle{czjphys209} 

\def\GLqtwo{{GL}_{q}(2)}
\def\GLhtwo{{GL}_{h}(2)}

\def\GLpqtwo{{GL}_{p,q}(2)}
\def\GLhh'two{{GL}_{h,h'}(2)}

\def\Grs{{G}_{r,s}}
\def\Gmk{{G}_{m,k}}

\def\rinv{{r}^{-1}}
\def\sinv{{s}^{-1}}
\def\Casmir{{\bf D}}
\def\Casmirinv{{\Casmir}^{-1}}
\def\ident{{\bf 1}}

\def\lmt{\longmapsto}

\def\cop{\Delta}

\def\cnt{\varepsilon}
\def\ot{\otimes}

\begin{document}
\title{REALISATIONS OF QUANTUM {\boldmath $\GLpqtwo$} AND JORDANIAN
{\boldmath $\GLhh'two$}}
\authori{DEEPAK PARASHAR and ROGER J. McDERMOTT}
\addressi{School of Computer and Mathematical Sciences, 
The Robert Gordon University,\\ 
St. Andrew Street, Aberdeen AB25 1HG, United Kingdom\\
(e-mail: deeps@scms.rgu.ac.uk , rm@scms.rgu.ac.uk)}
\authorii{}
\addressii{}
\authoriii{}  
\addressiii{} 
\headtitle{Realisations of Quantum $\GLpqtwo$ and
Jordanian $\GLhh'two$}
\headauthor{Deepak Parashar and Roger J. McDermott} 
\specialhead{Deepak Parashar: Realisations of Quantum $\GLpqtwo$ and
Jordanian $\GLhh'two$}

\evidence{A}
\daterec{XXX}    
\cislo{0}  \year{1999}
\setcounter{page}{1}
\pagesfromto{000--000}
\maketitle

\begin{abstract}
The quantum group $GL_{p,q}(2)$ is known to be related to the Jordanian
$GL_{h,h'}(2)$ via a contraction procedure. It can also be realised using
the generators of the Hopf algebra $G_{r,s}$. We contract the $G_{r,s}$
quantum group to obtain its Jordanian analogue $G_{m,k}$, which provides a
realisation of $GL_{h,h'}(2)$ in a manner similar to the $q$-deformed
case.
\end{abstract}

\section{Introduction}

Non Standard (or Jordanian) deformations of Lie groups and Lie algebras
has been a subject of considerable interest in the mathematical physics
community. Jordanian deformations for $GL(2)$ were introduced in [1,2],
its two parametric generalisation given in [3] and extended to the
supersymmetric case in [4]. Non Standard deformations of $sl(2)$ (i.e. at
the algebra level) were first proposed in [5], the universal
$R$-matrix presented in [6-8] and irreducible representations studied
in [9,10]. A peculiar feature of this deformation (also known as
$h$-deformation) is that the corresponding $R$-matrix is triangular.
It was shown in [11] that up to isomorphism, $\GLqtwo$ and $\GLhtwo$ are
the only possible distinct deformations (with central determinant) of the
group $GL(2)$. In [12], an interesting observation was made that the
$h$-deformation could be obtained by a singular limit of a similarity
transformation from the $q$-deformations of the group $GL(2)$. Given this
contraction procedure, it would be useful to look for Jordanian
deformations of other $q$-groups. 
\par

In the present paper, we focus our attention on a particular two parameter
quantum group, denoted $\Grs$, which provides a realisation of the well
known $\GLpqtwo$. We investigate the contraction procedure on $\Grs$, in
order to obtain its non standard counterpart. The generators of the
contracted structure are employed to realise the two parameter non
standard $\GLhh'two$. This is similar to what happens in the $q$-deformed
case. 
\par

\section{Quantum {\boldmath $\Grs$} and Realisation of {\boldmath 
$\GLpqtwo$}}

The two parameter quantum group $\Grs$ is generated by elements $a$, $b$,
$c$, $d$, and
$f$ satisfying the relations
\[
\begin{array}{ll}
ab=\rinv ba,\quad &db=rbd\\
ac=\rinv ca,\quad &dc=rcd\\
bc=cb,\quad &[a,d]=(\rinv-r)bc
\end{array}
\]
and
\[
\begin{array}{ll}
af=fa,\quad &cf=sfc\\
bf=\sinv fb,\quad &df=fd
\end{array}
\]
Elements $a$, $b$, $c$, $d$ satisfying the first set of
commutation relations form a subalgebra which coincides exactly with
$\GLqtwo$ when $q = \rinv$. The matrix of generators is 
\[
T=\left(\begin{array}{ccc}a&b&0\\c&d&0\\0&0&f\end{array}\right)
\]
and the Hopf structure is given as
\[
\begin{array}{lll}
\cop (T)&=& T\dot{\ot} T\\
\cnt (T)&=& \ident
\end{array}
\]
The Casimir operator is defined as $\Casmir = ad-\rinv bc$. The inverse is
assumed to exist and satisfies $\cop(\Casmirinv)=\Casmirinv\ot 
\Casmirinv$, $\cnt (\Casmirinv)=1$, $S(\Casmirinv)=\Casmir$, which enables
determination of the antipode matrix $S(T)$, as
\[
S\left(\begin{array}{ccc}a&b&0\\c&d&0\\0&0&f\end{array}\right)
=\Casmirinv
\left(\begin{array}{ccc}d&-rc&0\\-\rinv c&a&0\\0&0&\Casmir
f^{-1}\end{array}\right)
\]
The quantum determinant $\delta = \Casmir f$ is group-like but not
central.

The quantum group $\Grs$ was proposed in [13] as a particular quotient of
the multiparameter $q$-deformation of $GL(3)$. The structure of $\Grs$ is
interesting because it contains the one parameter $q$-deformation of
$GL(2)$ as a Hopf subalgebra and also gives a simple realisation of the
quantum group $\GLpqtwo$ in terms of the generators of $\Grs$. There is a
Hopf algebra morphism $\mathcal{F}$ from
$\Grs$ to $\GLpqtwo$ given by
\[
{\mathcal{F}}: \Grs\lmt \GLpqtwo
\]
\[
{\mathcal{F}}\left(\begin{array}{cc}a&b\\c&d\end{array}\right)\lmt
\left(\begin{array}{cc}a'&b'\\c'&d'\end{array}\right)=
f^{N}\left(\begin{array}{cc}a&b\\c&d\end{array}\right)
\]

The elements $a'$,$b'$,$c'$ and $d'$ are the generators of $\GLpqtwo$ and
$N$ is a fixed non-zero integer. The relation between the deformation
parameters $(p,q)$ and $(r,s)$ is given by
\[p = \rinv s^{N} \quad , \quad \quad q = \rinv s^{-N}\]

This quantum group can, therefore, be used to realise both $\GLqtwo$ and
$\GLpqtwo$ quantum groups.

\newpage
\section{{\boldmath $R$}-matrices and Contraction limits}

The $R$- matrix of $\Grs$ explicitly reads

\[
R(\Grs) = \left(\begin{array}{ccccccccc}
r & 0 & 0 & 0 & 0 & 0 & 0 & 0 & 0\\
0 & 1 & 0 & r-r^{-1} & 0 & 0 & 0 & 0 & 0\\
0 & 0 & s & 0 & 0 & 0 & r-r^{-1} & 0 & 0\\
0 & 0 & 0 & 1 & 0 & 0 & 0 & 0 & 0\\
0 & 0 & 0 & 0 & r & 0 & 0 & 0 & 0\\
0 & 0 & 0 & 0 & 0 & 1 & 0 & r-r^{-1} & 0\\
0 & 0 & 0 & 0 & 0 & 0 & s^{-1} & 0 & 0\\
0 & 0 & 0 & 0 & 0 & 0 & 0 & 1 & 0\\
0 & 0 & 0 & 0 & 0 & 0 & 0 & 0 & r
\end{array}\right)
\]

with entries labelled in the usual numerical order $(11)$, $(12)$, $(13)$,
$(21)$, $(22)$, $(23)$, $(31)$, $(32)$, $(33)$. If we reorder the indices
of this $R$-matrix with the elements in the order $(11)$, $(12)$, $(21)$,
$(22)$, $(13)$, $(23)$, $(31)$, $(32)$, $(33)$, then we obtain a block
matrix, say $R_{q}$ which is similar to the form of the $\GLqtwo$
$R$-matrix with the $q$ in the $R^{11}_{11}$ position itself replaced by
the $\GLqtwo$ $R$-matrix.

\[
R_{q} = \left(\begin{array}{cccc}
R(GL_{r}(2)) & 0 & 0 & 0\\
0 & S & \lambda I & 0\\
0 & 0 & S^{-1} & 0\\
0 & 0 & 0 & r
\end{array}\right)
\]

where $R(GL_{r}(2))$ is the $4\times 4$ $R$-matrix for $\GLqtwo$ with
$q=r$, $\lambda=r-r^{-1}$, $I$ is the $2\times 2$ identity matrix and $S$
is the $2\times 2$ matrix 
$S = \left( \begin{array}{cc}s&0\\0&1\end{array} \right)$ 
where $r$ and $s$ are the deformation parameters. The zeroes are the zero
matrices of appropriate order. The usual block structure of the $R$-matrix
is clearly visible in this form.

It is well known [12] that the non standard $R$-matrix $R_{h}(2)$
can be obtained from the $q$-deformed $R_{q}(2)$ as a singular limit of a
similarity transformation
\[
R_{h}(2) = \lim_{q\rightarrow 1}(g^{-1}\ot g^{-1}) R_{q}(2)(g\ot g)
\]
where $g=\left( \begin{array}{cc}1&\eta\\0&1\end{array} \right)$.
Such a transformation has been generalised to higher dimensions [14] and
has also been successfully applied to two parameter quantum groups.
Here we apply the above transformation for the $\Grs$ quantum group. Our
starting point is the block diagonal form of the $\Grs$ $R$-matrix,
denoted $R_{q}$
\[
R_{q} = \left(\begin{array}{ccccccccc}
r & 0 & 0 & 0 & 0 & 0 & 0 & 0 & 0\\
0 & 1 & \lambda & 0 & 0 & 0 & 0 & 0 & 0\\
0 & 0 & 1 & 0 & 0 & 0 & 0 & 0 & 0\\
0 & 0 & 0 & r & 0 & 0 & 0 & 0 & 0\\
0 & 0 & 0 & 0 & s & 0 & \lambda & 0 & 0\\
0 & 0 & 0 & 0 & 0 & 1 & 0 & \lambda & 0\\
0 & 0 & 0 & 0 & 0 & 0 & s^{-1} & 0 & 0\\
0 & 0 & 0 & 0 & 0 & 0 & 0 & 1 & 0\\
0 & 0 & 0 & 0 & 0 & 0 & 0 & 0 & r
\end{array}\right)
\]
where $\lambda=r-r^{-1}$. We apply to $R_{q}$ the transformation
\[
(G^{-1}\ot G^{-1})R_{q}(G\ot G)
\]
Here the transformation matrix $G$ is a $3\times 3$ matrix and chosen in
the block diagonal form
\[
G=\left(\begin{array}{cc}g&0\\0&1\end{array}\right)
\]
where $g$ is the transformation matrix for the two dimensional
case. Substituting $\eta=\frac{m}{r-1}$ and then taking the singular limit
$r\rightarrow 1$, $s\rightarrow 1$ (such that $\frac{1-s}{1-r}\rightarrow
\frac{k}{m}$) yields the Jordanian  $R$-matrix
\[
R_{h} = R(\Gmk) = \left(\begin{array}{ccccccccc}
1 & m & -m & m^{2} & 0 & 0 & 0 & 0 & 0\\
0 & 1 & 0 & m & 0 & 0 & 0 & 0 & 0\\
0 & 0 & 1 & -m & 0 & 0 & 0 & 0 & 0\\
0 & 0 & 0 & 1 & 0 & 0 & 0 & 0 & 0\\
0 & 0 & 0 & 0 & 1 & k & 0 & 0 & 0\\
0 & 0 & 0 & 0 & 0 & 1 & 0 & 0 & 0\\
0 & 0 & 0 & 0 & 0 & 0 & 1 & -k & 0\\
0 & 0 & 0 & 0 & 0 & 0 & 0 & 1 & 0\\
0 & 0 & 0 & 0 & 0 & 0 & 0 & 0 & 1
\end{array}\right)
\]
where the entries are labelled in the block diagonal form $(11)$, $(12)$,
$(21)$, $(22)$, $(13)$, $(23)$, $(31)$, $(32)$, $(33)$. It is
straightforward to verify that this $R$-matrix is triangular and a
solution of the Quantum Yang Baxter Equation
\[
R_{12}R_{13}R_{23} = R_{23}R_{13}R_{12}
\]
It is interesting to note that the block diagonal form of $R(\Gmk)$ embeds
the $R$-matrix for the single parameter deformed $\GLhtwo$ for $m=h$.

\section{Jordanian {\boldmath $\Gmk$} and Realisation of {\boldmath
$\GLhh'two$}}

A two parameter Jordanian quantum group, denoted $\Gmk$, can be formed by
using the contracted $R$-matrix $R(\Gmk)$ in conjunction with a $T$-matrix
of the form
\[
T=\left(\begin{array}{ccc}a&b&0\\c&d&0\\0&0&f\end{array}\right)
\]
The $RTT$- relations give the commutation
relations between the generators $a$, $b$, $c$, $d$ and $f$.
\[
\begin{array}{llll}
&[c,d] = mc^{2},& \qquad &[c,b] = m(ac+cd) = m(ca+dc)\\
&[c,a] = mc^{2},& \qquad &[d,a] = m(d-a)c = mc(d-a)
\end{array}
\]
\[
[d,b] = m(d^{2}-D)
\]
\[
[b,a] = m(D-a^{2})
\]
and
\[
\begin{array}{lll}
&[f,a] = kcf, \qquad &[f,b] = k(df-fa)\\
&[f,c] = 0, \qquad &[f,d] = -kcf
\end{array}
\]
The element $D = ad-bc-mac = ad-cb+mcd$ is central in the whole algebra. 
The coalgebra structure of $\Gmk$ can be written as
\[
\begin{array}{lll}
\cop (T)&=& T\dot{\ot} T\\
\cnt (T)&=& \ident
\end{array}
\]
Adjoining the element $D^{-1}$ to the algebra enables determination of the
antipode matrix $S(T)$,
\[
S\left(\begin{array}{ccc}a&b&0\\c&d&0\\0&0&f\end{array}\right)
=D^{-1}
\left(\begin{array}{ccc}
d-mc & -b-m(d-a)+m^{2}c & 0\\
-c & a+mc & 0\\
0&0&Df^{-1}\end{array}\right)
\]
(The Hopf structure of $D^{-1}$ is $\cop(D^{-1})=D^{-1}\ot D^{-1},
\cnt(D^{-1})=\ident, S(D^{-1})=D$.)

It is evident that the elements $a$, $b$, $c$ and $d$ of $\Gmk$ form a
Hopf subalgebra which coincides with non standard $GL(2)$ with deformation
parameter $m$. This is exactly analogous to the $q$-deformed case where
the first four elements of $\Grs$ form the $\GLqtwo$ Hopf subalgebra.
Again, the remaining fifth element $f$ generates the $GL(1)$ group, as it
did in the $q$-deformed case, and the second parameter appears only
through the cross commutation relations between $GL_{m}(2)$ and $GL(1)$
elements. Therefore, $\Gmk$ can also be considered as a two parameter
Jordanian deformation of classical $GL(2)\ot GL(1)$ group.
\par

Now we wish to explore the connection of $\Gmk$ with the two parameter
Jordanian $\GLhh'two$. A Hopf algebra morphism 
\[
{\mathcal{F}}: \Gmk\lmt \GLhh'two
\]
of exactly the same form as in the $q$-deformed case, exists between the
generators of $\Gmk$ and $\GLhh'two$ provided that the two sets of
deformation parameters $(h,h')$ and $(m,k)$ are related via the equation
\[
h = m + Nk \quad , \quad \quad h' = m - Nk
\]
Note that for vanishing $k$, one gets the one parameter case. In addition,
using the above realisation together with the coproduct, counit and
antipode axioms for the $\Gmk$ algebra and the respective homeomorphism
properties, one can easily recover the standard coproduct, counit and
antipode for $\GLhh'two$. Thus, the non standard $\GLhh'two$ group can in
fact be reproduced from the newly defined non standard $\Gmk$. It is
curious to note that if we write $p=e^{h}$, $q=e^{h'}$, $r=e^{-m}$ and
$s=e^{k}$, then the relations between the parameters in the $q$-deformed
case and the $h$-deformed case are identical. 
\par

\section{Conclusions}

We have applied the contraction procedure to the $\Grs$
quantum group and obtained a new Jordanian quantum group $\Gmk$. The
group $\Gmk$ has five generators and two deformation parameters and
contains the single parameter $\GLhtwo$ as a Hopf subalgebra. Furthermore,
we have given a realisation of the two parameter $\GLhh'two$ through the
generators of $\Gmk$ which also reproduces its full Hopf algebra
structure. The results match with the $q$-deformed case.
\par

\end{document}